\documentclass[letterpaper, 10 pt, conference]{ieeeconf}  

\IEEEoverridecommandlockouts                              

\usepackage{amsfonts}
\usepackage{tikz}
\usetikzlibrary{positioning, arrows.meta}
\usepackage{subcaption}
\usepackage{graphicx}
\usepackage{capt-of}

\newcommand{\edit}[1]{{\color{black}#1}}
\newcommand{\diag}{\mathbf{diag}}

\usepackage{amsmath} 
\usepackage{amssymb}  

\title{\LARGE \bf
Large-Scale Network Utility Maximization via GPU-Accelerated Proximal Message Passing}

\author{Akshay Sreekumar$^{1}$, Anthony Degleris$^{2}$, Ram Rajagopal$^{1}$
\thanks{*This work was not supported by any organization}
\thanks{$^{1}$Akshay Sreekumar and Ram Rajagopal are with the Department of Electrical Engineering,
        Stanford University, 350 Jane Stanford Way, Stanford CA, USA
        {\tt\small \{akshay81,ramr\}@stanford.edu}}%
\thanks{$^{2}$Anthony Degleris is with Gridmatic Inc, 20450 Stevens Creek Boulevard, Cupertino CA, USA
        {\tt\small anthony@gridmatic.com}}%
}
\begin{document}

\maketitle
\thispagestyle{empty}
\pagestyle{empty}

\bstctlcite{IEEEexample:BSTcontrol}

\begin{abstract}

We present a GPU-accelerated proximal message passing algorithm for large-scale network utility maximization (NUM). NUM is a fundamental problem in resource allocation, where resources are allocated across various streams in a network to maximize total utility while respecting link capacity constraints. Our method, a variant of ADMM, requires only sparse matrix–vector multiplies with the link–route matrix and element-wise proximal operator evaluations, enabling fully parallel updates across streams and links. It also supports heterogeneous utility types, including logarithmic utilities common in NUM, and does not assume strict concavity. We implement our method in PyTorch and demonstrate its performance on problems with tens of millions of variables and constraints, achieving 4$\times$ to 20$\times$ speedups over existing CPU and GPU solvers and solving problem sizes that exhaust the memory of baseline methods. We also show that our algorithm is robust to congestion and link-capacity degradation. Finally, a time-expanded transit seat allocation case study illustrates how our approach yields interpretable allocations in realistic networks.

\end{abstract}

\section{INTRODUCTION}

Network Utility Maximization (NUM) is an optimization framework for allocating resources across competing users in a network. Kelly et al. studied and formalized the connection between congestion control and fair resource allocation in large networks such as the Internet \cite{kelly_rate_1998}. This paradigm for fair resource allocation in networks is general and also has applications to sensor networks, caching systems, and transit systems~\cite{deng_maximizing_2016, dehghan_utility_2016, yin_maximizing_2022}. We focus on Internet and transit networks, which regularly scale to $10^5-10^6$ variables and constraints
\cite{shakkottai_network_2008}. 

Methods for solving NUM range from centralized interior point methods \cite{boyd_interior-point_2007} to distributed algorithms, many of which are based on dual decomposition \cite{wei_distributed_2013, palomar_tutorial_2006, bickson_distributed_2009}. These dual decomposition methods typically require \textit{strictly} concave utilities. The logarithmic utilities considered in this work allow NUM to be reformulated as a conic optimization problem with exponential cone constraints. State-of-the-art commercial solvers, e.g.\ MOSEK, support interior-point algorithms with native exponential cone handling, and are well suited for solving small to medium instances of NUM \cite{mosek_aps_mosek_2025}. Recently, GPU-accelerated algorithms have emerged for large-scale convex optimization. CuClarabel is a GPU-accelerated interior-point method for solving conic problems, capable of handling exponential cone constraints \cite{chen_cuclarabel_2024}. However, CuClarabel relies on \edit{cuDSS} to solve linear systems, which can be prohibitively slow for very large-scale problems. The authors in \cite{applegate_practical_2021} develop a variant of the primal-dual hybrid gradient (PDHG) algorithm to efficiently solve large-scale linear programs. Several works leverage GPU acceleration for restarted PDHG to solve linear and quadratic programs \cite{lu_cupdlpjl_2024, lu_mpax_2025}. While these PDHG methods do not involve solving a linear system, instead requiring only matrix-vector products, none of the previous solvers support the logarithmic objectives that are commonly used in NUM. The authors in \cite{zhang_solving_2025} develop a GPU-compatible algorithm based on PDHG for solving large multicommodity network flow problems. Kraning et al. develop the proximal message passing (PMP) algorithm, a variant of the alternating direction method of multipliers (ADMM), for solving very large DC optimal power flow problems (DC-OPF) on electricity networks~\cite{kraning_dynamic_2013}. Degleris et al. extend the framework of \cite{kraning_dynamic_2013} by implementing the PMP algorithm for DC-OPF on the GPU, using only sparse incidence matrix multiplies and vectorized scalar operations \cite{degleris_gpu_2024}. 

In this work, we adapt the GPU-accelerated PMP algorithm of \cite{degleris_gpu_2024} to efficiently solve large-scale NUM. Our method is fully distributed, handles extremely large problem instances, and can support multiple utility functions.

\section{PROBLEM SETTING}
We consider a NUM problem involving $n$ traffic streams and $m$ links. Each traffic stream $j$ has a fixed route comprised of some subset of the links and a utility function \edit{$U_j : \mathbb{R} \to \mathbb{R} \cup \{\infty \}$} that is \edit{the sum of a concave, twice differentiable function and an indicator function over a convex set.} The utility derived from assigning a stream rate $x_j$ is given by $U_j(x_j)$. In this work, we consider utility functions \edit{$U_j(x_j)$ with linear ($w_jx_j$) or logarithmic ($w_j\log(x_j$)) components.} The latter corresponds to \textit{weighted proportional fairness}, a special case of the $\alpha$-fair family of utilities studied in \cite{mo_fair_1998}. Our framework is amenable to other choices of utility functions as well, provided they have a simple to compute proximal operator. One such example is $\alpha=2$ fairness, which corresponds to \textit{minimum potential delay fairness} in communication networks \cite{shakkottai_network_2008}. 
The total network utility is $U(x) = \sum_{j=1}^n U_j(x_j)$,
where $x \in \mathbb{R}^n$ is the vector of stream rates.
Let $R \in \mathbb{R}^{m \times n}$ be the link-route matrix, defined as
\begin{equation}
R_{ij} =
\begin{cases}
1 & \text{if stream } j \text{ uses link } i, \\
0 & \text{otherwise}.
\end{cases}
\end{equation}
Each link $i \in \{1, \dots, m\}$ has capacity $c_i > 0$. Traffic on link $i$ must not exceed $c_i$, and stream rates must be nonnegative. 

The NUM problem is
\begin{equation}
\begin{array}{ll}
    \underset{\edit{x}}{\textrm{maximize}} \quad & \sum_{j=1}^n U_j(x_j) \\
    \text{subject to} \quad
    & Rx \leq c, \\
    & x \geq 0,
\end{array}
\label{eq:vanilla_num}
\end{equation}
where the variable is $x \in \mathbb{R}^n$.
We can re-write Prob.~\eqref{eq:vanilla_num} as
\begin{equation}
\begin{array}{ll}
    \underset{\edit{x,s}}{\textrm{minimize}} \quad & \sum_{j=1}^n -U_j(x_j) \\
    \text{subject to} \quad
    & Rx + s = c, \\
    & s \geq 0, \quad\edit{x \geq 0}
    
\end{array}
\label{eq:num_w_slack}
\end{equation}
where $s\in\mathbb{R}^m$ is a slack variable. When the utilities $U_j(x_j)$ are concave, Prob.~\eqref{eq:num_w_slack} is a convex minimization problem. 

\section{Proximal Message Passing}

To derive the PMP algorithm for NUM problems, we view Prob.~\eqref{eq:num_w_slack} as a bipartite graph where the two node sets are streams and links, and the edges that connect them are \textit{terminals}. Let the set of terminals be $\mathcal{J}=\{1,\ldots,J\}$, the set of links $\mathcal{L}=\{l_1,\ldots,l_{\edit{m}}\}$, and the set of streams $\mathcal{S}=\{\sigma_1,\ldots, \edit{\sigma_{n}, \ldots, \sigma_{m+n}}\}$, \edit{where the first $n$ streams are the traffic streams and the remaining $m$ streams are the slack streams.}
\edit{Here, $J = \mathrm{nnz}(R) + m$, where $\mathrm{nnz}(R)$ is the number of non-zeros in $R$.} 
Each terminal $j\in\mathcal J$ is incident to precisely one stream $\sigma\in\mathcal S$ and one link $l\in\mathcal L$, so the stream set and the link set each partition $\mathcal J$. 
We will show that the structure of Problem~\eqref{eq:num_w_slack} can be interpreted as this bipartite graph with streams, links, and terminals. Each terminal in the bipartite graph contains a \textit{flow} between streams and links. Streams are nodes with costs for producing or consuming flow, while links are nodes with flow-balance constraints.

To illustrate this, we expand a single row of equality constraints to understand what exactly are the flows that are summing to 0 at each link. Doing so for row $i$, we obtain
\begin{equation}
\begin{array}{l}
\sum_{j} R_{ij}x_{j} + s_{i} - c_{i} = 0.
\label{eq:expanded_net_row}
\end{array}
\end{equation}
Note that $x,s$ are the optimization variables in our problem. From~\eqref{eq:expanded_net_row} we can see two types of flow terms: $R_{ij}x_j$, which represent the flows of traffic stream $j$ connected to link $i$, and $s_i - c_i$, which can be interpreted as the flow of a ``slack'' stream at link $i$. We always have two forms of streams: (i) the original traffic streams and (ii) slack streams. Slack streams allow us to artificially saturate the capacity of a link such that flow balance holds. 

We define the variable $p \in \mathbb{R}^{|\mathcal{J}|}$ to denote the terminal flows in the system. We write $j \in \sigma$ and $j \in l$ to indicate that terminal $j$ is connected to stream $\sigma$ or link $l$. Let $|\sigma|$ and $|l|$ denote the number of terminals connected to a given stream or link. We use set-valued indices such as $p_l \in \mathbb{R}^{|l|}$ and $p_\sigma \in \mathbb{R}^{|\sigma|}$ to represent the flows associated with all terminals connected to link $l$ and stream $\sigma$, respectively.

\subsection{A Small Example}
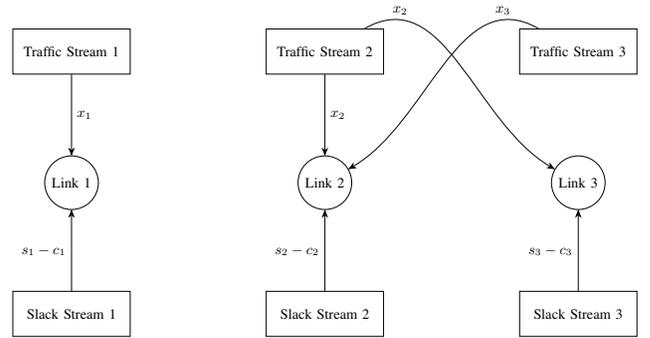
\begin{figure}[t]
\centering
\scalebox{0.6}{
\begin{tikzpicture}[
    >=Stealth,
    every node/.style={font=\small},
    device/.style={draw, minimum width=2.6cm, minimum height=1cm},
    slack/.style={draw, minimum width=2.6cm, minimum height=1cm},
    net/.style={draw, circle, minimum size=1cm},
    node distance=1.5cm
]

\node[device] (D1) {Traffic Stream 1};
\node[device, right=1.2cm of D1] (D2) {Traffic Stream 2};
\node[device, right=1.2cm of D2] (D3) {Traffic Stream 3};
\node[device, right=1.2cm of D3] (D4) {Traffic Stream 4};

\node[net, below=1.8cm of D1, xshift=1.65cm] (N1) {Link 1};
\node[net, below=1.8cm of D2, xshift=1.65cm] (N2) {Link 2};
\node[net, below=1.8cm of D3, xshift=1.65cm] (N3) {Link 3};

\node[slack, below=1.8cm of N1, xshift=0cm] (S1) {Slack Stream 1};
\node[slack, below=1.8cm of N2, xshift=0cm] (S2) {Slack Stream 2};
\node[slack, below=1.8cm of N3, xshift=0cm] (S3) {Slack Stream 3};

\draw[->] (S1) -- node[left] {$s_1 - c_1$}  (N1);
\draw[->] (S2) -- node[left] {$s_2 - c_2$}  (N2);
\draw[->] (S3) -- node[left] {$s_3 - c_3$}  (N3);

\draw[<-] (N1) -- node[right] {$x_1$} (D1);
\draw[<-] (N2) -- node[right] {$x_2$} (D2);

\draw[<-] (N2) to[out=0, in=270] 
    node[pos=0.9,  right] {$x_3$} 
    (D3);

\draw[<-] (N3.west) to[out=180, in=0] 
    node[pos=0.9, above right] {$x_2$} 
    (D2.east);

\draw[<-] (N3) to[out=0, in=270] 
    node[pos=0.9, left] {$x_4$} 
    (D4.south);
\end{tikzpicture}
}
\caption{Bipartite Graph for NUM}
\label{fig:concrete_network}
\end{figure}

Consider a network with \edit{four} traffic streams and three links with the following set of equality constraints:
\begin{equation*}
\begin{array}{ll}
\begin{pmatrix}
1 & 0      & 0 & 0 \\
0      & 1 & 1 & 0 \\
0      & 1 & 0 & 1
\end{pmatrix}
\begin{pmatrix}
x_1 \\ x_2 \\ x_3 \\ x_4
\end{pmatrix}
+
\begin{pmatrix}
s_1 \\ s_2 \\ s_3
\end{pmatrix}
-
\begin{pmatrix}
c_1 \\ c_2 \\ c_3
\end{pmatrix}
=
\begin{pmatrix}
0 \\ 0 \\ 0
\end{pmatrix}
\label{eq:concrete_equality_constrs}
\end{array}
\end{equation*}
Fig. \ref{fig:concrete_network} shows the corresponding bipartite graph. 
\subsection{Streams}
Every stream (either a \emph{traffic} or a \emph{slack} stream) has a utility function. \edit{We define the convex indicator function $\mathcal{I}_\mathcal{C}(\cdot)$ as}

\begin{equation*}
    \edit
    {\mathcal{I}_{\mathcal{C}}(x) :=}
\begin{cases}
\edit
{0}, & \edit{ \text{if } \quad x  \in \mathcal{C},}\\
\edit{+\infty,} & \edit{\text{otherwise},}
\end{cases}
\end{equation*}

\edit{where $\mathcal{C}$ is a convex set.}

We now define the \edit{utility} functions for various streams.


\paragraph{Log-Utility Stream}

\begin{equation}
\label{eq:log_utility_device_cost_func}
\begin{array}{ll}
f_{\sigma}(p_{\sigma})
= \displaystyle{\min_{\edit{x_\sigma \geq 0}}} \left(
\; -\,w_\sigma \log x_\sigma
  + \edit{\mathcal{I}_{\{0\}}\{\,p_\sigma - \mathbf{1}_{|\sigma|} x_\sigma\,\}} \right)
\end{array}
\end{equation}
Here, $\edit{x_\sigma \in \mathbb{R}_{+}}$ acts as a local variable for the traffic stream \edit{flow} that controls the terminal flows. The constraint $p_{\sigma}=\mathbf{1}_{|\sigma|} x_\sigma$, \edit{where $\mathbf{1}_{|\sigma|}$ is the all-ones vector of length $|\sigma|$}, means that for a traffic stream with multiple terminals, the flow in each incident terminal must exactly equal $x_\sigma$ \edit{(i.e. all entries of $p_{\sigma}$ are identical)}. The scalar $w_\sigma > 0$ is the stream weight.

\paragraph{Linear-Utility Stream}

\begin{equation}
\label{eq:linear_utility_device_cost_func}
\begin{array}{ll}
f_{\sigma}(p_{\sigma})
= \displaystyle{\min_{\edit{x_\sigma \geq 0}}} \left(
\; -\,w_\sigma x_\sigma
  + \edit{\mathcal{I}_{\{0\}}\{\,p_\sigma - \mathbf{1}_{|\sigma|} x_\sigma\,\}} \right)
\end{array}
\end{equation}

\paragraph{Slack Stream}

\begin{equation}
\begin{array}{ll}
        f_{\sigma}(p_{\sigma}) = \edit{\mathcal{I}_{\mathbb{R}_+}\{p_{\sigma} + c_{\sigma} \}}
    \end{array}
\label{eq:slack_device_cost_func}
\end{equation}

\paragraph{Other Utility Streams}
Our method also accommodates more complex utility functions. We only require \edit{these utility functions to have} simple proximal operators (see Section~\ref{sec:proximal_updates}).

\subsection{Links}

The sum of all flows at a link must balance to 0. \edit{For every link $l\in\mathcal L$, we define the function $g_{l}(\cdot)$ as an indicator function to capture this flow balance constraint. That is, we have $g_l(p_l) := \mathcal{I}_{\{0\}}(\mathbf{1}^\top p_{l})$.}

\subsection{Proximal Message Passing from ADMM}

Prob.~\eqref{eq:num_w_slack} can be reformulated as
\begin{equation}
\begin{array}{ll}
    \underset{\edit{p,z}}{\textrm{minimize}} \quad & \sum_{\sigma \in \mathcal{S}} f_\sigma(p_\sigma) + \sum_{l \in \mathcal{L}} g_l(z_l) \\
    \text{subject to} \quad
    & p=z. \\
\end{array}
\label{eq:admm_form}
\end{equation}
We interpret Prob.~\eqref{eq:admm_form} as follows: $n$ utility streams generate flows on their route links, and $m$ slack streams inject the residual flow needed to saturate each link’s capacity. The first objective term encodes stream costs while the second enforces feasibility via link flow balance constraints. 
\edit{We can recover the solution to~\eqref{eq:num_w_slack} via the relations $x_\sigma \mathbf{1}_{|\sigma|} =  p_\sigma$ for $\sigma \leq n$ and $s_{\sigma - n} = p_\sigma + c_{\sigma - n}$ for $n < \sigma \leq n+m$.}

As is typical in ADMM, we introduce a copy variable $z$. The variable $p$ represents a stream-side copy of the terminal flow variable, while $z$ is the link-side copy of the flow variables. Our subproblems in ADMM are separable across streams and links. Following \cite{kraning_dynamic_2013, degleris_gpu_2024, boyd_distributed_2010} we can write the scaled augmented Lagrangian for Prob.~\eqref{eq:admm_form} as
\begin{equation}
\textstyle L_\rho(p, z, u) = \sum_{\sigma \in \mathcal{S}} f_\sigma(p_\sigma) + \sum_{l \in \mathcal{L}} g_l(z_l) + \frac{\rho}{2} \left\|p - z + u\right\|_2^2,
\label{eq:augmented_lagrangian}
\end{equation}

where $\rho > 0$ is a penalty parameter, and $u$ is the scaled dual variable given by $y/\rho$ where $y$ is the unscaled dual variable. \edit{The scaled dual variables for the streams and links are $u_{\sigma}$ and $u_{l}$ respectively.} Using the augmented Lagrangian in~\eqref{eq:augmented_lagrangian}, we can derive the ADMM iterations as
\begin{equation*}
\textstyle
 p_\sigma^{k+1} := \displaystyle{\mathop{\arg\min}_{p_\sigma}}
 \textstyle
 \left( f_\sigma(p_\sigma) + \frac{\rho}{2} \left\| p_\sigma - z_\sigma^k + u_\sigma^k \right\|_2^2 \right), \; \sigma \in \mathcal{S}
\end{equation*}

\begin{equation*}
\textstyle
z_l^{k+1} := \displaystyle{\mathop{\arg\min}_{z_l}} 
\textstyle
\left( g_l(z_l) + \frac{\rho}{2} \left\| z_l - u_l^k - p_l^{k+1} \right\|_2^2 \right), \; l \in \mathcal{L}
\end{equation*}
\begin{equation*}
\textstyle
u_l^{k+1} := u_l^k + \left( p_l^{k+1} - z_l^{k+1} \right), \; l \in \mathcal{L}.
\end{equation*}

\edit{By solving a constrained least squares problem, we can show that $z_l^{k+1} = u_l^{k} + \bar p_l^{k+1}$ and rewrite the above updates as the \textit{proximal message passing algorithm}}.
\begin{enumerate}
    \item \textit{Proximal Device Updates}
    \begin{equation}
    \begin{array}{ll}
        p_\sigma^{k+1} := \mathrm{prox}_{f_\sigma, \rho} \left( p_\sigma^k - \bar{p}_\sigma^k - u_\sigma^k \right), \; \sigma \in \mathcal{S}
        \label{eq:prox_updates}
        \end{array}
    \end{equation}
    \item \textit{Scaled Price Updates}
    \begin{equation}
    \begin{array}{ll}
        u_l^{k+1} := u_l^k + \bar{p}_l^{k+1}, \; l \in \mathcal{L}
        \end{array}
    \end{equation}
\end{enumerate}
where the proximal operator for a function $f$ is given by
\begin{equation}
\begin{array}{ll}
    \mathrm{prox}_{f, \rho}(\zeta) = \displaystyle{\mathop{\arg\min}_{y}} \left( f(y) + \frac{\rho}{2} \left\| y - \zeta \right\|_2^2 \right),
    \end{array}
\end{equation}
and $\bar{p}_\sigma \in \mathbb{R}^{|\sigma|}$ and $\bar{p}_l \in \mathbb{R}^{|l|}$ are vectors where each value is the average of all terminal flows incident to stream $\sigma$ or link $l$, respectively.
\edit{See \cite{kraning_dynamic_2013} for a detailed derivation.}

\subsection{Proximal Updates}
\label{sec:proximal_updates}

The efficacy of the proximal updates in \eqref{eq:prox_updates} depends on how efficiently we can compute \edit{$\mathrm{prox}_{f_\sigma, \rho}$} for the various streams in the system.
Fortunately, for traffic streams with logarithmic or linear utilities, as well as slack streams, \edit{the minimization problem in the proximal operator are one dimensional and can be solved analytically.}

\begin{enumerate}
    \item \textit{Log-utility streams:}
\begin{equation}
\begin{array}{ll}
\edit{\mathrm{prox}_{f_\sigma, \rho}(\zeta)} = \frac{\mathbf{1}_{|\sigma|}^T\zeta + \sqrt{(\mathbf{1}_{|\sigma|}^T\zeta)^2 + 4w_\sigma\ |\sigma| / \rho}} {2|\sigma|} \mathbf{1}_\sigma
\end{array}
\end{equation}

\item \textit{Linear-utility streams:}
\begin{equation}
\begin{array}{ll}
    \edit{\mathrm{prox}_{f_\sigma, \rho}(\zeta)} = \frac{\mathbf{1}_{|\sigma|}^{T}\zeta +\frac{w_\sigma}{\rho}}{|\sigma|} \mathbf{1}_\sigma
\end{array}
\label{eq:variable_device_optimal_soln}
\end{equation}

\item \textit{Slack streams:}
\begin{equation}
\begin{array}{ll}
    \edit{\mathrm{prox}_{f_\sigma, \rho}(\zeta)} = \max\{\zeta, -c_\sigma\}
\end{array}
\end{equation}
\end{enumerate}

\subsection{Convergence}

The primal and dual residuals at iteration $i$ are $r^{(i)} = \bar{p}^{(i)}$ and $s^{(i)} = \rho ((p^{(i)} - \bar{p}^{(i)}) - (p^{(i-1)} - \bar{p}^{(i-1)}))$, respectively. The primal residual is the net flow imbalance across all the links, which is exactly the measure of primal infeasibility in Prob.~\eqref{eq:admm_form}. The dual residual is the difference between consecutive iterates of the difference between the flows and average flow on each net. We terminate the PMP algorithm in accordance with the following criterion:
\begin{equation}
\begin{array}{ll}
\|r^{(i)}\|_2 < \varepsilon_{\text{tol}}, \quad \|s^{(i)}\|_2 < \varepsilon_{\text{tol}},
\end{array}
\end{equation}
where $\varepsilon_{\text{tol}} = \varepsilon_{\text{abs}} \sqrt{|\mathcal{J}|}$ and $\varepsilon_{\text{abs}} > 0$ is an absolute tolerance. Because PMP is a variant of ADMM, we refer the reader to \cite{boyd_distributed_2010} for additional details on ADMM's convergence. 

\subsection{Accelerating Convergence}

\paragraph{Residual Balancing}

We adaptively adjust the penalty parameter $\rho$ so as to balance the primal and dual residuals. We utilize the simple update scheme from \cite{boyd_distributed_2010}:
\begin{equation}
\textstyle
\rho^{(i+1)} :=
\begin{cases}
\gamma \rho^{(i)} & \text{if } \|r^{(i)}\|_2 > \mu \|s^{(i)}\|_2, \\[4pt]
\rho^{(i)}/\gamma & \text{if } \|s^{(i)}\|_2 > \mu \|r^{(i)}\|_2, \\[4pt]
\rho^{(i)} & \text{otherwise},
\end{cases}
\label{eq:rho-update}
\end{equation}
where $\gamma > 1$ and $\mu > 1$ are parameters of the update rule.

Each time we update $\rho$, we accordingly rescale $u$ as 
\begin{equation}
\begin{array}{ll}
u^{(i+1)} := \left( \frac{\rho^{(i)}}{\rho^{(i+1)}} \right) u^{(i+1)}.
\label{eq:dual-rescale}
\end{array}
\end{equation}
We set $\mu=2$ and $\gamma=1.1$ and update $\rho$ every $50$ iterations.

\paragraph{Over-relaxation}
Over-relaxation can improve the convergence of ADMM~\cite{boyd_distributed_2010}. Let $\alpha \in [1,2]$ denote the \emph{relaxation parameter}. In the $z,u$ updates, we replace $p^{k+1}_l$ with
\begin{equation}
\begin{array}{ll}
p^{k+1}_{l,+} \;:=\; \alpha\, p^{k+1}_l + (1-\alpha)\, z^{k}_l.
\end{array}
\end{equation}
With over-relaxation, our updates no longer simplify to the message passing updates derived above. Instead, we write
\begin{equation}
\begin{aligned}
\textstyle
z^{k+1}_l &= \alpha\big(p^{k+1}_l - \bar{p}^{k+1}_l\big) + (1-\alpha)\,z^{k}_l,\\
u^{k+1}_l &= u^{k}_l + \alpha\,\bar{p}^{k+1}_l.
\end{aligned}
\end{equation}
In our experiments, we find that $\alpha=1.6$ works well.

\begin{figure*}[t]
  \centering
  \begin{minipage}[t]{0.565\textwidth}
    \centering
    \includegraphics[width=\linewidth]{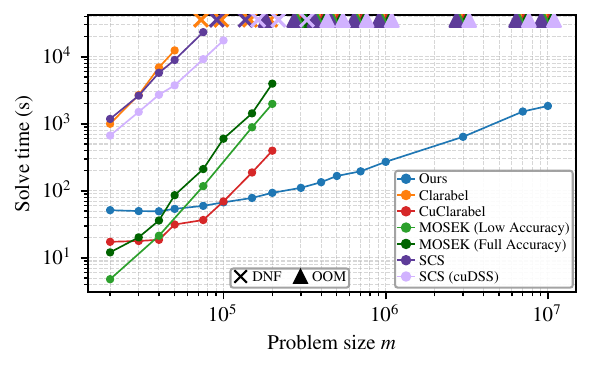}
    \captionof{figure}{Scaling performance of solvers on various problem sizes}
    \label{fig:solver_scaling}
  \end{minipage}\hfill
  \begin{minipage}[t]{0.425\textwidth}
    \centering
    \includegraphics[width=\linewidth]{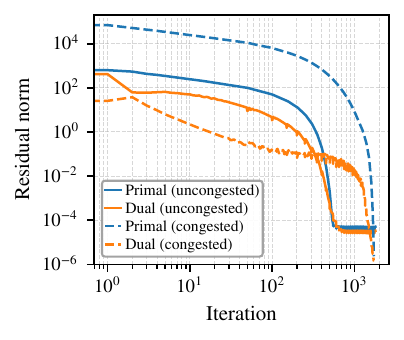}
    \captionof{figure}{Uncongested vs. Congested Convergence}
    \label{fig:congested_convergence}
  \end{minipage}
\end{figure*}

\section{GPU IMPLEMENTATION}
The two key operations in PMP that require efficient GPU implementation are (i) proximal updates for each of the streams and (ii) computing averages across links.

\subsection{Stream Data Model}
Streams are grouped into types $c$ by utility function and terminal count. This allows for vectorized proximal updates, since streams of the same type share the same utility function and dimensionality of flows. For each type $c$, we store its variables \edit{in a tensor of size $\tau_c \times |c|$}, where $|c|$ is the number of streams of that type \edit{and $\tau_c$ is the number of terminals per stream for that type}.

\subsection{Vectorization}

The proximal updates for the streams are closed-form, meaning we can compute them with fully vectorized, parallel updates. For each stream type $c$, we batch its $|c|$ streams and update them simultaneously. For example, the linear-utility update given in \eqref{eq:variable_device_optimal_soln} for a batch of $|c|=k$ streams is
\begin{equation}
\begin{array}{ll}
x^{*}_{\tau_c}
=
\Big(
\diag(\mathbf{J}_{\tau_c,k}^{\top} Z)
-
\diag\!\big(\tfrac{1}{\rho} w_{\tau_c}\big)
\Big)\,
\frac{1}{\tau_c} I_k,
\label{eq:matrix_prox_variable_device}
\end{array}
\end{equation}
where $\mathbf{J}_{\tau_c,k}\!\in\!\mathbb{R}^{\tau_c\times k}$ is the all-ones matrix, 
$Z=[z_1,\dots,z_k]\!\in\!\mathbb{R}^{\tau_c\times k}$ the stacked evaluation points (one column per stream), 
and $w_{\tau_c}\!\in\!\mathbb{R}^k$ the stream weights. 
Eq.~\eqref{eq:matrix_prox_variable_device} admits an efficient GPU implementation via broadcasting.

\subsection{Averages}

The PMP updates require efficiently computing the averages of all terminals incident to each link. These averages are identical for all terminals connected to the same link, so we store them as tensors of size~$m$. Conceptually, averaging gathers information across links and can be expressed as multiplication by a sparse incidence matrix, which can be implemented on the GPU using a scatter kernel.

Formally, for stream type $c$, let $R_{c_i} \in \mathbb{R}^{m \times |c|}$ denote the incidence matrix for terminal $i$, with entries
\begin{equation*}
(R_{c_i})_{l\sigma} =
\begin{cases}
1 & \text{if terminal $i$ of stream $\sigma$ connects to link $l$}, \\
0 & \text{otherwise},
\end{cases}
\end{equation*}
where $i = 1, \ldots, \tau_c$. The matrices $R_{c_i}$ are the incidence matrices for each terminal of each stream type. The average quantity $\bar{p}$ is then given by
\begin{equation*}
\begin{array}{ll}
\bar{p}
=
\frac{1}{|l|} \odot
\sum_{c} \sum_{i=1}^{\tau_c}
R_{c_i} \, p_{c_i},
\end{array}
\end{equation*}
where $p_{c_i}$ is a $|c|$-dimensional tensor, $\bar{p}$ is an $m$-dimensional tensor, $|l| \in \mathbb{R}^m$ is the vector of terminal counts per link, and $\odot$ denotes elementwise multiplication.

\subsection{Software Implementation}
We develop an open-source PyTorch implementation\footnote{https://github.com/degleris1/zap} of the PMP algorithm for NUM problems. Our implementation is modular, supports both linear and logarithmic utilities, and can be extended to other stream types by specifying a desired utility function and its corresponding proximal operator. On a single GPU, we can solve a problem with \edit{$m=1\times10^7$} links and \edit{$n=5\times10^6$} streams in $1847$ seconds. 

\begin{figure*}[t]
  \centering
  \includegraphics[width=\textwidth]{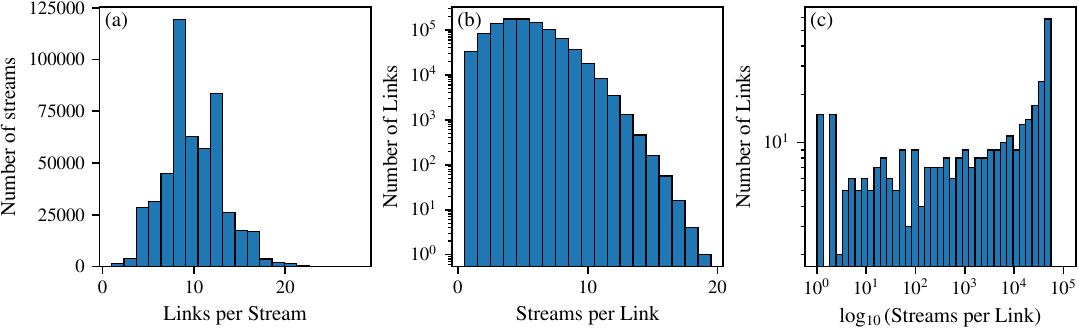}
  \caption{(a) Uncongested route lengths. 
           (b) Uncongested link utilization. 
           (c) Congested link utilization.}
  \label{fig:route_link_hists}
\end{figure*}

\begin{figure}[!t]
  \centering

    \includegraphics[width=\linewidth]{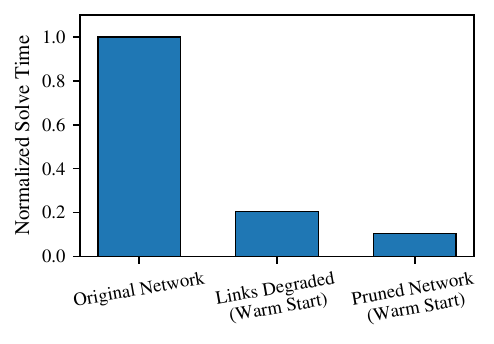}
    \caption{Warm starting under link degradation and failure}
    \label{fig:warm_starts}
    \vspace{-6pt} 
\end{figure}

\section{PERFORMANCE RESULTS}
We perform several numerical experiments to demonstrate the efficacy of our method on a variety of different problem sizes and scenarios. All experiments are performed on a single NVIDIA A100-SXM4 GPU with 80GB of memory, supported by 32 virtual CPU cores and 64GB of RAM. While
our experiments are synthetic in construction, they map to applications in Internet congestion
control and transit seat allocation. We first benchmark our method against several open-source and commercial solvers to demonstrate the superior scaling performance of our method on uncongested networks. We examine the performance of PMP in more complicated networks with congested links, showing that our method is robust and still able to converge to acceptable solutions. Finally, we observe that warm-starting our method allows it to quickly re-compute optimal allocations in the event of network degradation.

\subsection{Applications to Internet and Transit Settings}
The numerical experiments in this section are synthetic, but directly applicable to two practical domains. In Internet congestion control, streams are end-to-end sessions, links are bandwidth-limited communication links, and link degradation models outages. The goal is to maximize aggregate utility subject to bandwidth constraints. In transit, streams are itinerary–departure options and links are time-expanded vehicle–seat resources; congestion arises when itineraries overlap, and degradation reflects service disruptions.

\subsection{Scaling}
We randomly generate link–route matrices $R$ following the methodology of \cite{boyd_interior-point_2007, bickson_distributed_2009}. 
Given $m$ total links and $n=m/2$ traffic streams, we assign links to streams such that the average stream is comprised of 10 links. We benchmark our method against the commercial solver MOSEK and the open-source conic solvers \edit{CuClarabel, Clarabel, and SCS \cite{mosek_aps_mosek_2025, goulart_clarabel_2024, chen_cuclarabel_2024, odonoghue_operator_2021}. CuClarabel is designed for GPU acceleration; for SCS, we report results with both its standard linear system solver and also the cuDSS backend.} All solvers are run with low-accuracy settings unless otherwise stated, \edit{meaning we relax each solvers stopping criteria (e.g. primal/dual residuals, optimality gap) so that runs terminate at coarse accuracy.} Fig.~\ref{fig:solver_scaling} shows wall-clock time versus problem size, i.e. the number of links $m$. If a solver fails to return a solution within 10~hours, we mark it as Did Not Finish (DNF); if it fails due to memory limits, we mark it as Out of Memory (OOM). For small problem sizes, MOSEK and CuClarabel outperform our method, but at $m=2\times 10^5$ our approach is $\sim$4$\times$ faster than CuClarabel and $\sim$20$\times$ faster than MOSEK. For $m \geq 5\times 10^5$, only our method solves the instances without running out of memory. We also note that MOSEK on low accuracy fails on three instances ($m=3\times10^4$, $m=5\times 10^4$, and $m=10^5$), returning NaN.

\subsection{Congested Networks}
We model congestion by making a small fraction of links heavily utilized. We randomly select $0.1\%$ of links and connect each to approximately $10\%$ of streams. In Fig.~\ref{fig:congested_convergence}, we see PMP reaches medium accuracy with primal and dual residuals $<10^{-4}$ in about $1{,}000$ iterations in the uncongested case. Congestion reduces the sparsity of $R$ and slows convergence slightly. However, similar accuracy is achieved with only a few hundred additional iterations. 

In Fig.~\ref{fig:route_link_hists}a, we show the distribution of links per stream for an uncongested network with \edit{$m=1\times10^6$} and \edit{$n=5\times10^5$}. In Fig.~\ref{fig:route_link_hists}b–\ref{fig:route_link_hists}c, we then plot the distribution of streams per link for both uncongested and congested networks of the same dimension. Relative to the uncongested case, congestion introduces a heavy tail in the streams per link distribution. These heavily used links correspond to dense rows in $R$, causing the modest slowdown in convergence.

\subsection{Reallocation under Link Degradation}

Next, we consider optimal allocations under capacity \emph{degradation} of the links. Starting from a solution $x_0^\star$ for link-route matrix $R$ and capacities $c$, we re-solve under degraded capacities $c_d$ where each link independently has a $25\%$ chance of a $50\%$ capacity reduction. In Fig.~\ref{fig:warm_starts}, we see warm-starting from $x_0^\star$ yields a roughly $5\times$ speedup. 

In the extreme, links fail entirely. Given the structure of NUM, a failed link removes its row of $R$ and eliminates all streams that traverse it. We then solve the pruned instance  with $R_{\text{pruned}}, c_{\text{pruned}}$ after failing each link with probability $25\%$. In Fig.~\ref{fig:warm_starts}, we see combining pruning and warm-starting yields a roughly $10\times$ speedup in solve time. 

\section{CASE STUDY: TRANSIT SEAT ALLOCATION}
\label{sec:seat_allocation}

\subsection{Problem Setting}
We consider a trip–reservation setting \cite{yin_maximizing_2022} with multiple itinerary–departure options. We model a transit system as a time-expanded network with stations $\mathcal{V}$, directed spatial edges $\mathcal{E}\subseteq\mathcal{V}\times\mathcal{V}$, and discrete time bins $t=0,\dots,T\!-\!1$ of length $\Delta t$. A \emph{link} is a time-stamped edge, $\mathcal{L}=\mathcal{E}\times\{0,\dots,T\!-\!1\}$ with $\ell=(e,t)$. The travel time for all links is $\tau_\ell=\Delta t$. 

An itinerary–departure \emph{stream} is $j=(k,r,t_0)$, where $(o_k,d_k)$ is an OD pair, $r\in R_k$ is a spatial route ($R_k$ is the set of all spatial routes for OD pair $k$), and $t_0$ is a departure bin. Each stream receives allocation $x_j > 0$ and utility $U_j(x_j)=w_j\log x_j$, with $w_j$ encoding itinerary preferences (e.g., off-peak travel). Links represent scheduled vehicle–seat capacity and streams represent passenger options. We solve the continuous relaxation and round as needed.

\begin{figure}[t]
  \centering
  \subfloat[Low and high streams]{%
    \includegraphics[width=0.44\columnwidth]{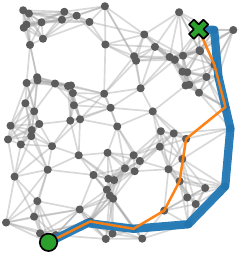}%
    \label{fig:flows}}
  \hfill
  \subfloat[Normalized link price $\hat{\lambda}$]{%
    \includegraphics[width=0.54\columnwidth]{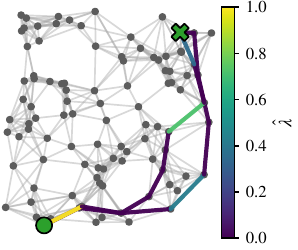}%
    \label{fig:prices}}
  \caption{Flows and link prices for two routes departing at $t_0=0$.}
  \label{fig:routes-two-small}
\end{figure}

\subsection{Instance Generation and Solution}

We build a synthetic time-expanded network with $S=100$ stations and time step $\Delta t=5$ minutes, with $T=192$. This yields $m=|E|\,T=182{,}784$ links and $n=19{,}800$ streams with weighted logarithmic utilities ($w_j=1$). To interpret a slice of the solution, we compare two itineraries for the \emph{same} OD and departure. Fig.~\ref{fig:flows} shows the spatial routes. Line thickness is proportional to the allocated flow. Here, the blue route receives nearly $5\times$ the flow.

The disparity follows from the link shadow prices $\lambda_{(\ell,t)}\!\ge 0$. We define the \emph{path price} of stream $j$ as: 
\begin{equation}
\begin{array}{ll}
\pi_j \;=\; \sum_{(\ell,t)} R_{(\ell,t),j}\,\lambda_{(\ell,t)} .
\end{array}
\end{equation}
For weighted log utilities, KKT stationarity gives $\frac{w_j}{x_j} \;=\; \pi_j$, so lower-priced paths receive more flow. In Fig.~\ref{fig:prices}, normalized prices along each route show an additional bottleneck on the orange path that raises $\pi_j$, whereas the blue path only has the same initial high price link. NUM concentrates price on congested links and shifts allocations to cheaper routes. 

\section{CONCLUSION}
We have presented a GPU-accelerated proximal message passing solver for large-scale NUM. Our method uses only sparse matrix-vector multiplies with the link–route matrix and closed-form proximal updates. We implement this method in pure PyTorch and show that it scales to millions of variables, delivers $4\times-20\times$ speedups over strong CPU/GPU baselines, and solves instances those methods cannot. We then demonstrate that our method is robust under congestion and link-capacity degradation and benefits significantly from warm starting after perturbations. Finally, we use a transit seat–allocation case study to demonstrate how NUM can be applied to a realistic operations problem, yielding interpretable allocations across overlapping itinerary options and underscoring our approach's broader applicability.

\bibliographystyle{IEEEtran}
\bibliography{ieee_ctrl, references}

\end{document}